\documentclass[12pt]{elsarticle}

\usepackage[compatible]{algpseudocode}
\usepackage{algorithm}

\floatname{algorithm}{Procedure}

 \usepackage{pgfplots}
 \usepackage{tikz}
 \usetikzlibrary{shapes,arrows}
  \usetikzlibrary{external}
 \tikzset{external/force remake}
\usepackage{amssymb}
\usepackage{makeidx}
\usepackage{amsmath}
\usepackage{graphicx}
\usepackage{amsfonts}
\usepackage{eurosym}
\usepackage{xifthen}
\usepackage{xcolor}
\newboolean{pedro}
\setboolean{pedro}{true}
\newcommand{\pedro}{\ifthenelse{\boolean{pedro}}{\color{black}
    \setboolean{pedro}{false}}{\color{black}\setboolean{pedro}{true}}}
\newboolean{luis}
\setboolean{luis}{true}
\newcommand{\luis}{\ifthenelse{\boolean{luis}}{\color{black}
    \setboolean{luis}{false}}{\color{black}\setboolean{luis}{true}}}

\title{Mid-term bio-economic optimization of multi-species fisheries}
\author{L. Bay\'on\corref{cor1}}
\author{P. Fortuny Ayuso\corref{cor2}}
\author{P. J. Garc\'ia-Nieto\corref{cor2}}
\author{J. A. Otero\corref{cor2}}
\author{P. M. Su\'arez\corref{cor2}}
\author{C. Tasis\corref{cor2}}

\address{Department of Mathematics, University of Oviedo, E.P.I. Campus of Viesques, Gij\'{o}n, 33203, Spain.}
\cortext[cor1]{Corresponding author: bayon@uniovi.es}

\newproof{pf}{Proof}

\journal{Applied Mathematical Modelling}
\begin{document}

\begin{abstract}
In this paper, we analyze the dynamics of a multi-species fisheries system in
the presence of harvesting. We solve the problem of finding the optimal
harvesting strategy for a mid-term horizon with a fixed final stock of
each species, while maximizing the expected present value of
total revenues. The problem is formulated as an optimal control problem. For
its solution, we combine techniques derived from Pontryagin's Maximum Principle,
cyclic coordinate descent and the shooting method. The algorithm we develop can
solve problems both with inter-species competition and with
predator-prey behaviors. Several numerical examples are presented to
illustrate the different possibilities of the method and a study of the dependence of the
behavior on some parameters is performed.
\end{abstract}


\maketitle

\section{Introduction}

A large number of problems in applied mathematics is concerned with the analysis
of dynamical processes in renewable resources. \luis{} A classical introductory
reference to the topic is the book \cite{Perman 2003}, whose chapters contain
abundant bibliography related to each of the traditionally considered stocks:
fish, forests, or freshwater. More recently, the excellent handbooks
\cite{Kennedy 2012} or \cite{Cardenas 2013} have appeared. Works centered on
each of the stocks are \cite{Beverton 2012} or \cite{Hilborn 2013} for the
dynamics of fish populations; \cite{Bettinger 2016} ---using Dynamic
Programming--- for the problem of forest management is studied and, quite
recently, \cite{Mays 2018} uses mathematical programming and differential
dynamic programming techniques  to study the optimization of groundwater
management systems, or freshwater inflows among other topics.\luis{}

\luis{}From the economic point of view, the renewable resources problem is
stated as the maximization of profit over a certain time horizon, where the
problem is subject to the biological dynamic of the resource, an initial stock
size of fish, and other technological or political constraints \cite{Sundstrom
2010}. The economic value includes, by its very nature, a
discount factor. The models related to fishery and
marine economics and its relation to the economic foundation are discussed in
\cite{Kronbak 2014}. In \cite{Grune 2015}, the authors present a new
approach to solve dynamic decision models in economics. Their is based
on an iterative solution of optimal control problems in finite time horizons.\luis{}

In this paper, we deal with the particular case of harvesting of fish
\cite{Prouzet 2015}. Its optimization in different ways is a complex problem
with both ecological and economical implications. The interaction between these
two aspects is nowadays of great interest, and a large corpus has been
written. \luis{}In \cite{Schluter 2012}, the so-called social-ecological
systems, characterized by strong links between both aspects, are analyzed. An analysis of how the technical tools should be used as a
basis for policy making and management within the European context is presented
in \cite{de Jonge 2012}. Also, \cite{Brown 2013} describes how the global
economy can be restructured to make it compatible with the Earth's
ecosystem. The study presented in \cite{Kasperski 2015} maximizes the net
present value of a multi-species fishery where species interact both
ecologically (in the ecosystem) and economically (in the output markets).  The
interactions taking place in an integrated ecological--economical model are also
analyzed in \cite{Speers 2016}.  \luis{}

From the biological point of view, most papers have focused on the optimal
harvesting policies considering only one fish species, and an infinite-time
horizon $[0,\infty)$ (see, for example, \cite{Poudel 2012}, \cite{Shah 2013}
and \cite{Qiu 2013}). These works deal mostly with long-term
stability and the existence of a steady-state solution. A
method developed by Tsur \cite{Tsur 2001} calculates the optimal steady-state,
and does not require the solution of the dynamic optimization problem.

Tsur's method was applied in \cite{Suri 2008} for a single-state,
infinite-horizon, autonomous (except for the discount factor) model. The same
problem was studied by some of the authors in a previous paper
\cite{Bayon 2017}
but including, besides the classical infinite-time horizon problem, an
example with long-term horizon. In the latter, we considered a fixed and finite
optimization interval $[0,T]$, and we assumed it to be long enough
for the steady state to be reached in its duration. We observed in those works
that the
solution leaves the steady state when approaching the end of the interval and
seeks a final stock value.
This optimal value is not part of the statement: it is freely sought
by the system.

\luis{}In recent years, the development of a proper model incorporating
multi-species fisheries has gained interest. The biological
interactions between species play an important role in optimal fisheries
management, and it will be one of the main aspects studied in this paper. This
multi-species modelization leads to complex dynamical problems. The most classical standpoint can be found, for instance, in \cite{Matsuda 2006} where the authors assume that the per capita growth
rate of each species is a linear function of population abundance. A more
sophisticated multi-species bioeconomic model is analyzed in \cite{Kellner
2011}. This model (of an archetypal Caribbean coral reef community)
includes two focal prey species (parrotfish and snapper) and a generalist
predator (grouper). Again, the predation mortality of the two preys is
modeled as linear on the predator population.

That type of linear model, with a food chain consisting of only two trophic
levels and a prey--predator system is frequently used in the literature (see \cite{Legovic 2010}). The same authors present, in \cite{Legovic 2012},
a system of $n$ cooperating populations without harvesting. In this
system of logistic populations, each one, by its mere presence, increases the
carrying capacity of the other ones: one species may be a
commensal to any number of species and likewise, one species may receive help
from any number of species.  \luis{}

To cope with these complexities, authors sometimes analyze mainly the
stability of the biological system in lack of fishing \cite{Kar 2011}. Other
times, they use specific models which are only valid for particular
groups of species (see \cite{Agnarsson 2008}, \cite{Aanestad 2009} and
\cite{Poudel 2015}), or they develop more general methods but always under the
infinite-horizon assumption and searching the optimal steady-state
\cite{Agmour 2017}.

\luis{} There is, then, a large variety of models. In the review paper
\cite{Plaganyi 2014}, an interesting classification of three classes of
ecological models used for an ecosystem approach to fisheries is presented:
Single-species assessments, Models of Intermediate Complexity for Ecosystem
assessments and Whole-of-ecosystem models. In the same reference, a subset of
realistic models is shown (with the appropriate references), with the Region and
the Species modelled. Sixteen other examples of population models with
explicitly accounted-for predation mortality are also given in \cite{Tyrrell
  2011}.

We end this brief survey citing \cite{Kissling 2012}, where the principal ways
of modelling interactions between multiple species are presented. These include
(among others): Competition, Predator-Prey, Host-Parasitoid, Infectious
Diseases, Plant Competition, Facilitation. These models aim to incorporate
multispecies interactions at large spatial extents using interaction
matrices.\luis{} Our approach differs radically: we propose first of all, a
rather general modelization for studying multi-species interactions, with which
we can deal both with pure competence between species and with predator-prey
systems (and mixed competence-predator-prey ones).

\luis{}
Another important issue which can be taken into account in complex decision problems is the influence of uncertainty over fisheries. We shall deal, in this work, with a deterministic environment but we cannot omit the citation of several studies which include stochasticity. The book \cite{Kennedy 2012}, includes a part dealing with fisheries management from the deterministic point of view and a comparison with the stochastic one. An analogue comparative study (but on forestry management) can be found in the same reference. The effects of uncertainty over fisheries are also covered in \cite{Hilborn 2013}. The uncertainties and climate change impacts on forests has been recently analyzed in \cite{Linder 2014}, and \cite{Yousefpour 2012} is an excellent review on forest management. Finally, \cite{Kvamsdal 2016} analyzes a continuous, nonlinear bioeconomic model to demonstrate how stochasticity in the growth of fish stocks affects the optimal exploitation policy when prices are stochastic.
\luis{}

Once the model, under a deterministic environment, is presented, we focus on an
original setting: a multi-species biological system (say with $n$ species) and a
mid-term finite horizon $[0,T]$. \pedro{}In this context, we
understand mid-term as a relatively short time in comparison with the time taken
by the system to reach, by itself, a steady state, and at the same time, long
enough that days or even a single fishing campaign are small in
comparison. Following our computations, this gives a value for $T$ of around
$10$ years\pedro{}. In this horizon, a single manager has a target for the end
of the optimization interval: the final stock of each species $x_i(T)$
($i=1,\ldots, n$) is fixed \emph{a priori}.  This problem \pedro{}of considering
a specific final stock\pedro{}, which has barely appeared in the literature, is
essential in today's society: rather than seeking long-term biological
equilibrium by leaving the system to itself, it is important for authorities to
have a tool enabling them to lead the system to a desired state, respecting its
biological laws, while at the same time maximizing the profits of fishing during
that time interval.

To solve the problem, we adapt an algorithm that
takes advantage of the main
result of Optimal Control: Pontryagin's Maximum Principle \cite{Chiang 2000},
combined with the shooting method \cite{Aronna 2013} and an adapted version
for functionals of the cyclic coordinate descent \cite{Luo 1992}. This
resolution method can deal with arbitrarily complex multidimensional problems.
The idea lies in not trying to solve the problem globally but to deal with
it as an iterative sequence of one-dimensional problems in which, at each
step, each species considers the others as known functions. This property,
together with other specifics which we shall show in detail, is what allows the method not to be gravely influenced by the dimensionality
of the problem or by
the complexity of the model.

The paper is organized as follows: Section 2 presents the
different modelings needed for the problem being addressed; in Section 3, the
optimization algorithm is shown;
Section 4 shows a numerical example to illustrate the
algorithm's performance and the dependence of the model on several
variables is illustrated; finally, Section 5 includes the conclusions and perspectives for further research.

\section{Statement of the problem}

\subsection{Biological model}
\pedro{}Denoting by $t$ the time variable (in units of years), we let $x(t)$ represent the stock of a single species at time $t$. We work in the continuous setting because we are modeling species which have large stocks (cod, capelin, herring, etc.) and because the fishing events, which are daily, happen at a very small scale relatively to our time interval (and the influence of a single event on the stock is very small). In this context, $x(t)$ will be a continuous function with piecewise continuous derivative (i.e. the plot of $x(t)$ may have corner points but no discontinuities).\pedro{}

The model for the pattern of biological growth of one resource. The literature
\cite{Perman 2003} presents several models in which the actual growth rate
\pedro{}with respect to time\pedro{} $\dot{x}(t)$, depends on the stock size
$x(t)$. A commonly used functional form is the simple logistic function or
Verhulst equation \eqref{Log1}, where \pedro{}$r>0$\pedro{} denotes
the intrinsic growth rate and \pedro{}$k>0$\pedro{} the carrying
capacity of the species. This model is a good approximation to the natural
growth processes of many fish populations and will be the one we use
henceforward. Thus, from now on $f_{l}(x)$ will denote the right hand side of
\eqref{Log1}:
\begin{equation}
\dot{x}(t)  =f_{l}(x)=rx(t)\left(  1-\frac{x(t)}{k}\right) \label{Log1}
\end{equation}
In this model, the maximum amount of growth or maximum sustainable yield,
$x_{MSY}$, happens when the stock size is equal to $k/2$.
We consider upper and lower limits for the stock $x(t)$, expressed as
$x\in\lbrack x_{\min},k]$. These are imposed, respectively, by the
biological minimum of the species, $x_{\min}$ which allows its reproduction
and by the value $k$ of the model.

There are other generalizations of the logistic growth model such as the
modified logistic model \eqref{Log2} or the Gompertz equation \eqref{Log3}. \pedro{}Model \eqref{Log2} is used, for example, in \cite{Agnarsson 2008} for cod in Norway, with $\gamma=2$\pedro{}
\begin{align}
\dot{x}(t)  & =rx^{\gamma}(t)\left(  1-\frac{x(t)}{k}\right)\pedro{};\;\gamma>1\pedro{} \label{Log2}\\
\dot{x}(t)  & =rx(t)\ln\frac{k}{x(t)}\label{Log3}%
\end{align}
There are also discrete-population models such as the Beverton-Holt model
(which can be considered as a discrete version of the simple logistic function),
the Hassell model or the Ricker model (see a review in \cite{Brannstrom 2005}).
We do not study any of these.

\subsection{Harvesting model}

Human harvesting is included in the model as follows: denote
by $h(t)$ the rate of biomass harvest,\pedro{}which we assume is a piecewise continuous function\pedro{}. For the sake of simplicity, many
authors (see, for example, \cite{Perman 2003}, \cite{Suri 2008} and \cite{Shah
2013}) assume that all the different dimensions of the harvesting activity
(e.g. size of nets, number of trawlers, number of fishing days) can be
aggregated into one single magnitude called effort, $E(t)$. Thus, a
first model proposed by them is:
\begin{equation}
h(t)=qE(t)x(t)
\end{equation}
where $q$ is a constant, often called the catchability coefficient,
$E(t)$ is the fishing effort, and $x(t)$ is the fish stock level at time $t$.
The proportionality constant $q$ encodes the \emph{easiness} of fish harvesting.
This approach has obvious advantages in terms of mathematical tractability
but, in our view, it constitutes an unnecessary simplification.


In this paper, we model the dynamics in a more general form as:%
\begin{equation}
\dot{x}(t)=f_{l}(x(t))-h(t)\label{St Eq}%
\end{equation}
where $h(t)$ does not depend on any other quantities.

\subsection{Economic Model}

The economic model for each species follows. Let $\pi(x(t),h(t))$
be the instantaneous net revenue from the harvest of the stock biomass, given
as in \cite{Sandal 1997}:
\begin{equation}
\pi(x(t),h(t))=p(h(t))h(t)-c(x(t),h(t))\label{Rev}%
\end{equation}
where $p(h(t))$ is the price function and $c(x(t),h(t))$ the cost
function associated
with the harvest. We assume an economic model verifying the following three
natural conditions:
\begin{equation}
  \frac{\partial p(h)}{\partial h}<0;
  \frac{\partial c(x,h)}{\partial h}>0;\text{
}\frac{\partial c(x,h)}{\partial x}<0
\end{equation}
With these assumptions, we adopt the very general models:
\begin{align}
p(h(t))  & =p_{0}-p_{1}h\label{p}\\
c(x(t),h(t))  & =\frac{ch(t)^{\alpha}}{x(t)}\label{c}%
\end{align}
where $p_{0}$ is the stock price, $p_{1}$ is the strength of demand,
$c$ is the cost of exploitation and $\alpha$ is the harvest cost parameter.
These models correspond to \pedro{}real-world fisheries\pedro{} (see \emph{\cite{Agnarsson
2008} \pedro{}where it is used for Norway, Iceland and Denmark, and\pedro{} \cite{Poudel 2015}\pedro{} where it is used for the Barents Sea\pedro{})}, where the price of the harvest depends on the
amount harvested and the cost of harvesting depends on the stock biomass.
Substituting (\ref{p}) and (\ref{c}) in (\ref{Rev}), the profit function for
each species is:%
\begin{equation}
\pi(x(t),h(t))=p_{0}h(t)-p_{1}h(t)^{2}-\frac{ch(t)^{\alpha}}{x(t)}\label{Rev2}%
\end{equation}
We remark that other authors use simpler models, in which either the stock price
is constant \cite{Agmour 2017}, or the cost is linear in the harvesting
($\alpha=1$) or even constant \cite{Kar 2011}. In \cite{Poudel 2015}, some
species have
constant price and the price of others depend on the harvesting.
As regards the cost of harvesting, in those
works it does not depend on the stock for some species. Our model, inspired
by that of \cite{Poudel 2012}, is the most complete.


\subsection{Multi-species biological model}
Finally, we present the multi-species model for the biological system.
\luis{}As stated in the Introduction, we consider the classical modeling (see,
for example \cite{Matsuda 2006}) as that in which the per capita growth rates of
each species is a linear function of population abundances:
\begin{equation}
\frac{dn_{i}}{dt}=\left( \sum_{j=1}^{s}a_{ij}n_{j}\right) n_{i}
\end{equation}%
where $a_{ij}$ is the effect of species $i$
on the per capita growth rate of species $j$. The value of $a_{ij}$
is constant in the model, negative if either species $j$ exploits species $i$ or species $j$ and $i$ compete directly with each other. Its value is positive if, on the contrary, species $i$ exploits species $j$ and $a_{ji}\leq -a_{ij}$.
One can also define, following \cite{Matsuda 2006} the \emph{energy conversion rate} from prey $i$ to predator $j$ as:
\begin{equation}
m_{ji}=-\frac{a_{ji}}{a_{ij}}
\end{equation}%

The pattern above is the most usual. See, for example \cite{Kar 2011} where two
mutually competing fish species are considered, both subject to harvesting, and
both of them have a common unharvested predator. A modified version is shown in the same reference, where a predator-prey system subject to
harvesting is considered using a modified version of the Leslie-Gower scheme: it
is assumed that the reduction in a predator population has a reciprocal relation
with per capita availability of its preferred food.

See \cite{Agnarsson 2008} for other more particular models in which linear functions for the relation between two species is also considered but with power models for each of them separately. In \cite{Poudel 2015}, three fish species (capelin, cod and herring) are studied and the pairwise interactions are (at times) considered as square roots. The problem with this models is that they are too specific and unfit for generalization.

In this work, we are going to adopt a model which is even more general than the
one used in \cite{Agmour 2017}, where three species interact in a linear way
but only pairwise: we referred to this type of relation schema in the
Introduction, when citing \cite{Kissling 2012}. In it, three kinds of biological
interactions are distinguished: (1) simple qualitative linkages between species;
(2) quantitative interaction coefficients that reflect the power of these
interactions; and (3) parameters which measure ---by way of a proxy environmental variable--- how the triple interactions affect each species.
All the approaches above can be described using interaction matrices. However, there is no study including all the possible parameters we are going to include.


We propose a model with the following parameters:
\begin{itemize}
\item The effect $c_{ij}$ of species $j$ on the per-capita growth rate of species $i$.
\item The simultaneous effect $c_{ijk}$ of both species $j$ and $k$ on the per-capita growth rate of species $i$.
\item The convexity parameter $\gamma_{ij}$ of the influence of species $j$ on species $i$. It needs not be true that effects are linear on the population: this parameter measures this (possible) non-linearity.
\item The convexity parameter $\gamma_{ijk}$ of the pairwise influences of species $j$ and $k$ on species $i$, measured in species $j$.
\item The \emph{catalysis level} $\beta_{ij}$ of the species $i$ when interacted by species $j$: it might happen that \emph{the rate at which $i$ incresases is lower than the rate at which the influence of $j$ does}.
\item The \emph{catalysis level} $\beta_{ijk}$ of species $i$ when interacted by both species $j,k$.
\end{itemize}

One should take into account that a linear model assumes explicitly that interactions between individuals are straightforward (one more individual means exactly the same change in the quantity of the interaction). This assumption, although highly successful for basic modeling is a very elementary description of real-world systems. Introducing non-linearity allows us to provide a richer structure. A very simple case might be that, after some point, the increase in a prey species is irrelevant for the increase in a predator one, at least when the size of the predator species is small (or vice versa). This is, for instance, what $\beta_{ij}$ and/or $\gamma_{ij}$ can model, taking their values as $>1$ or $<1$. This is just for pairwise interactions. For triads (or larger) of interactions, our model copes perfectly well.
\luis{}
Hence, in the general case of dimension $n$, we model the dynamics of the $i$-th
species as:%
\begin{equation}
\dot{x}_{i}(t)=f_{l,i}(x_{i})-g_{i}(t,x_{1},\ldots ,x_{i},\ldots
,x_{n})-h_{i}(t)  \label{DYN}
\end{equation}

where $\mathbf{x}(t)=(x_{1}(t),...,x_{n}(t))$ and $g_{i}(\mathbf{x}(t))$ is
the function representing the competence or \textquotedblleft
coupling\textquotedblright\ among the different species, which we assume is
of the form:%
\begin{equation}
g_{i}(t,x_{1},\ldots ,x_{i},\ldots ,x_{n})=\sum_{\substack{ 1\leq j\leq n \\ %
j\neq i}}c_{ij}x_{i}(t)^{\beta _{ij}}x_{j}(t)^{\gamma _{ij}}+\sum_{\substack{
1\leq j<k\leq n \\ j\neq i\neq k}}c_{ijk}x_{i}(t)^{\beta
_{ijk}}x_{j}(t)^{\gamma _{ijk}}x_{k}(t)^{\gamma _{ikj}}  \label{G}
\end{equation}

This function we are proposing expresses the relations between those species
which compete with the $i$-th one. It is the most general one we have found
in the literature. It allows modeling with a high degree of detail,
including, for instance, simultaneous interactions of up to $3$ species.
Notice that this value might be increased without a substantial modification
of the complexity of the underlying problem. Each species $x_{i}$ is
influenced by its exponent $\beta _{ij}$ when affected by species $x_{j}$.


If $\beta _{ij}=\gamma _{ij}=1,$ $\beta _{ijk}=\gamma _{ijk}=0$ one obtains
the classical models \cite{Agnarsson 2008}, \cite{Kar 2011}, \cite{Agmour
2017}. In \cite{Poudel 2015} one also finds simultaneous interactions
between $3$ species, as in our model (but his is a specific case, with the
constants already fixed). Notice that (\ref{G}) can deal both with problems
of pure competence between species and with predator-prey models. In the
former ones, the coefficients $c_{ij}$ and $c_{ijk}$ are positive, whereas
in the latter ones, predators have negative $c_{ij}$ and $c_{ijk}$ and preys
positive.\bigskip

In order to clarify the notation used in \eqref{G}, we show now, by
way of example, the ground problem we shall analyze later on and which
corresponds to the three-species case, $n=3$ which, as seen in \cite{Kar
2011}, already constitutes a highly dynamic system. Using the simple logistic
function for $f_{l,i}$, we get:

\begin{multline}
\dot{x}_{1}(t)=r_{1}x_{1}(t)\left( 1-\dfrac{x_{1}(t)}{k_{1}}\right) -
\label{SD1} \\
-[c_{12}x_{1}(t)^{\beta _{12}}x_{2}(t)^{\gamma _{12}}+c_{13}x_{1}(t)^{\beta
_{13}}x_{3}(t)^{\gamma _{13}}+c_{123}x_{1}(t)^{\beta _{123}}x_{2}(t)^{\gamma
_{123}}x_{3}(t)^{\gamma _{132}}]-h_{1}(t)
\end{multline}%
\begin{multline}
\dot{x}_{2}(t)=r_{2}x_{2}(t)\left( 1-\dfrac{x_{2}(t)}{k_{2}}\right) -
\label{SD2} \\
-[c_{21}x_{2}(t)^{\beta _{21}}x_{1}(t)^{\gamma _{21}}+c_{23}x_{2}(t)^{\beta
_{23}}x_{3}(t)^{\gamma _{23}}+c_{213}x_{2}(t)^{\beta _{213}}x_{1}(t)^{\gamma
_{213}}x_{3}(t)^{\gamma _{231}}]-h_{2}(t)
\end{multline}%
\begin{multline}
\dot{x}_{3}(t)=r_{3}x_{3}(t)\left( 1-\dfrac{x_{3}(t)}{k_{3}}\right) -
\label{SD3} \\
-[c_{31}x_{3}(t)^{\beta _{31}}x_{1}(t)^{\gamma _{31}}+c_{32}x_{3}(t)^{\beta
_{32}}x_{2}(t)^{\gamma _{32}}+c_{312}x_{3}(t)^{\beta _{312}}x_{1}(t)^{\gamma
_{312}}x_{2}(t)^{\gamma _{321}}]-h_{3}(t)
\end{multline}

Once this complex modelization is adopted, we remark that the choice of the
method we are going to develop in the next section is not just casual. Our
resolution method is not intrinsically affected by the complexity of the
model: as we use a cyclic coordinate descent method, what we do is to solve
a sequence of one-dimensional problems. In each of these, each species
``considers'' the other ones as fixed, so
that the modeling of the coupling term is, in this case, of the simple form
$g_{i}(t,x_{i}).$

\subsection{Objective Functional}

Our model considers an open-access fishery model, in which a single manager
takes the market price of fish as given. The manager's objective is to
maximize profits from the harvest schedule of the multi-species over a finite
time horizon $[0,T]$, fixing a target for the end of the optimization
interval: the final stock of each species, $x_{i}(T)$. Moreover,
the solution is subject to the dynamic constraint equations \eqref{DYN} and
other natural and policy restrictions involving limits for the harvest and
the stock. Hence, \pedro{}letting $\pi_i(x(t),h(t))$ as in Equation \eqref{Rev2} be the profit function for each species $i$\pedro{}, our problem is:
\begin{equation}
\begin{array}
[c]{ll}
& \underset{h_{i}(t)}{\max}%
{\displaystyle\int_{0}^{T}}
{\displaystyle\sum\limits_{i=1}^{n}}
\pi_{i}(x_{i}(t),h_{i}(t))e^{-\delta t}dt\\
\text{s.t.} & \dot{x}_{i}(t)=f_{l,i}(x_{i}(t))-g_{i}(x_{1}(t),\ldots
,x_{i}(t),\ldots,x_{n}(t))-h_{i}(t)\\
& x_{i}(0)=x_{i0};\text{ }x_{i}(T)=x_{iT}\\
& h_{i}(t)\in H_{i}=[h_{i\min},h_{i\max}]\\
& i=1,...,n
\end{array}
\label{PR}%
\end{equation}
where $\delta>0$ is the discount rate, $x_{i0}$ is the initial stock level and
$x_{iT}$ the final desired stock level.

\section{Optimization Algorithm}

Notice that Problem \eqref{PR} can be
stated as a Multi-dimensional Optimal Control
Problem. The numerical algorithm we propose for solving it
uses a particular strategy related to the cyclic
coordinate descent (CCD) method \cite{Luo 1992}. This method
minimizes a function of $n$ variables cyclically with respect to the
coordinates. With our method, the Multi-dimensional problem with $n$ species
can be solved as a sequence of problems with one single species assuming the
others fixed (what we shall call the Unidimensional Problem). Thus,
starting from an admissible solution of the multi-dimensional
problem, $\mathbf{s}^{0}$, we
compute a sequence $(\mathbf{s}^{k})$ and the algorithm will
calculate:%
\begin{equation}
\underset{k\rightarrow\infty}{\lim}\mathbf{s}^{k}%
\end{equation}

We shall also give a necessary condition of maximum for the Unidimensional
Problem from which we shall compute the solution with the following
two steps:
the construction of $x_{i}^{K}$ for $i=1,...,n$ and the calculation of the optimal
$K$.

With this method, the problem can be solved like a sequence of problems
each of whose error functional converges to zero.

\subsection{Step 1: The Multi-dimensional Problem}

As explained above, Problem \eqref{PR} can be stated as a Multi-dimensional
Optimal Control Problem simply setting the stock as the state variable,
$\mathbf{x}(t)=(x_{1}(t),\ldots,x_{n}(t))$ and the harvesting
$\mathbf{u}(t)=(h_{1}(t),\ldots,h_{n}(t))$ as the control variable. Taking into
account that $\dot{x}_{i}(t)=f_{l,i}(x_{i})-g_{i}(\mathbf{x})-h_{i}(t)$, we can
write%
\begin{equation}
h_{i}(t)=f_{l,i}(x_{i})-g_{i}(\mathbf{x(t)})-\dot{x}_{i}(t)=G_{i}(t,\mathbf{x}%
(t),\dot{x}_{i}(t))
\end{equation}
and impose the specified constraints on both variables.

\pedro{}
In what follows, a function $F(t,\mathbf{x},\mathbf{u})$ (the integrand) will be given. Whenever a path $\mathbf{x}(t)=\mathbf{s}(t)$ and a coordinate $i$ are fixed, we shall denote $F_i^{\mathbf{s}}$ the corresponding function $F_i^{\mathbf{s}}(t)=F(t,\mathbf{s}(t),\mathbf{u}(t))$. There is also a vector function $\mathbf{f}(t)=(f_1(t),\dots,f_n(t))$.

\emph{Remark 1}: The following hypotheses are assumed:
(i) $F_{i}^{\mathbf{s}}$ and $f_{i}$ are continuous; (ii) $F_{i}^{\mathbf{s}}$
and $f_{i}$ have continuous partial first derivatives with respect to $t$ and
$x_{i}$; (iii) the control $u_{i}(t)$ is piecewise continuous; (iv) the state
variable $x_{i}(t)$ is continuous and its derivative is piecewise continuous
(i.e. $x_{i}(t)$ admits corner points). The set of admissible controls is
compact and convex. These hypothesis arise from the natural assumptions of our
continuous model. As is customary, the set of differentiable functions of a real
variable in $[0,T]$ with continuous derivative is denoted $\hat{C}^1[0,T]$.
\pedro{}

The objective is to maximize the discounted profits given the
state dynamics and fixing a final value for the stock:
\begin{equation}
\underset{\mathbf{u}(t)}{\max}J=\int_{0}^{T}F\left(  t,\mathbf{x}%
(t),\mathbf{u}(t)\right)  dt\label{Funct}%
\end{equation}
\pedro{}where $\mathbf{u}(t)=(h_1(t),\dots,h_n(t))$,
subject to satisfying, for each $t\in[0,T]$\pedro{}:
\begin{align}
&\mathbf{\dot{x}}(t)  =\pedro{}\mathbf{f}\pedro{}\left(  t,\mathbf{x}(t),\mathbf{u}(t)\right)  ;\text{
 }\mathbf{x}(0)=\mathbf{x}_{0};\text{  }\mathbf{x}(T)=\mathbf{x}%
                      _{T}\label{state eq}\\
  &\pedro{}h_{1\min}\leq h_1(t)\leq h_{1\max},\pedro{}
  \pedro{}\dots,h_{n\min}\leq h_n(t)\leq h_{n\max}\pedro{}
\end{align}
Hence, we say that $\mathbf{x}(t)=(x_{1}(t),\ldots,x_{n}(t))$ is
an admissible solution if it belongs to the set
\begin{equation}
\Theta:=\{\mathbf{z}\in\left(  \hat{C}^{1}[0,T]\right)  ^{n} /\text{ }%
\mathbf{z}(0)=\mathbf{x}_{0}\text{, }\mathbf{z}(T)=\mathbf{x}_{T},\text{
}h_{i\min}\leq
G_{i}(t,\mathbf{z}(t),\dot{z}_{i}(t))
\leq h_{i\max},\text{
}i=1,\ldots,n\}
\end{equation}

This problem presents the following remarkable features:
First, it is a multidimensional problem. Second, the optimization interval is
finite and the final stock is fixed. Third, the time $t$ is not explicitly
present in the problem (it is a time-autonomous problem), except in the discount
factor and finally ---fourth--- constraints on the control are imposed. Due to
the nature of problem, Optimal Control Theory, and more specifically
Pontryagin's Maximum Principle \cite{Chiang 2000} applies straightforwardly.

As we explained in the Introduction, we are not interested, in this work,
in the equilibrium or steady-state solution, at which
the resource stock size is unchanging over time (a biological equilibrium) and
the harvesting is constant. We aim to find the dynamic solution,
i.e. the adjustment path towards the final desired state
starting from a given initial value.

From the  Multi-dimensional problem with $n$ species we construct a
sequence of problems with one single species, assuming the others fixed.

Let $\mathbf{s=}(s_{1},\cdots,s_{n})\mathbf{\in}$ $\Theta$, with $ u_{i}%
(t)=h_{i}(t)=G_{i}(t,\mathbf{s}(t),\dot{s}_{i}(t))$ and set:
\begin{equation}
F_{i}^{\mathbf{s}}(t,x_{i},u_{i}):=F(s_{1}(t),\cdots,s_{i-1}(t),x_{i}%
(t),s_{i+1}(t),\cdots,s_{n}(t),h_{1}(t),\cdots,h_{i-1}(t),u_{i}(t),h_{i+1}%
(t),\cdots,h_{n}(t))
\end{equation}
The problem which consists in finding:
\begin{equation}
\underset{u_{i}(t)}{\max}J_{i}^{\mathbf{s}}(x_{i}):=\int_{0}^{T}%
F_{i}^{\mathbf{s}}(t,x_{i},u_{i})dt\label{UFunct}%
\end{equation}

in the set
\begin{multline}
\Theta_{i}^{\mathbf{s}}:=\{z_{i}\in \pedro{}\hat{C}\pedro{}^{1}[0,T]/\text{ }z_{i}(0)=x_{i0}\text{,
}z_{i}(T)=x_{iT},\text{ }\label{Ucon}\\
h_{i\min}\leq G_{i}(s_{1}(t),\cdots,s_{i-1}(t),x_{i}(t),s_{i+1}(t),\cdots
,s_{n}(t),\dot{x}_{i}(t))\leq h_{i\max}\}
\end{multline}

is (certainly) a unidimensional problem.

We advanced in the Introduction that the solution to the Multi-dimensional
problem will be obtained by applying the following algorithm:

\pedro{}
\textbf{Remark 2.}
Notice that Problem \eqref{UFunct} with the conditions \eqref{Ucon} is a unidimensional optimization problem which admits a single solution $s^{\ast}_i(t)$.
\pedro{}

\pedro{}
\textbf{Definition 1.}
\emph{We define the $i$-th maximizing map $\varphi_i$ in $\Theta$ as the map which substitutes the $i$-th state function by the corresponding solution to problem \eqref{UFunct}, leaving the others unmodified: let $(s_1, \dots, s_n)\in \Theta$ be an admissible element. Let $s_i^{\ast}$ be the solution of Problem \eqref{UFunct} subject to \eqref{Ucon}. Then}
\begin{equation}
\varphi_{i}(s_{1},\dots,s_{i},\dots,s_{n})=(s_{1},\dots,s_{i}^{\ast}%
,\dots,s_{n}),
\end{equation}
\emph{so that $J_i^{\mathbf{s}}(s^{\ast})\geq J_i^{\mathbf{s}}(x_i)$ for all $x_i\in \Theta^{\mathbf{s}}_i$ (i.e. $s^{\ast}_i$ maximizes $J^{\mathbf{s}}_i$).}
\pedro{}


We shall denote by $\varphi$ the map associated with the ascent algorithm, which
will be the composition of the i-th maximizing map:
\begin{equation}
\varphi:=\varphi_{n}\circ\cdots\circ\varphi_{1}%
\end{equation}
In every $k$-th iteration of the algorithm, ``the $n$ components of
$\mathbf{s}$ will have been maximized'' by means of the $i$-th maximizing
applications for each $i=1,\ldots, n$ in this order,
thus obtaining the new admissible element $\mathbf{s}^{k}$:
\begin{equation}
\mathbf{s}^{k}=\varphi(\mathbf{s}^{k-1})=(\varphi_{_{n}}\circ\varphi_{_{n-1}}%
\circ\cdots\circ\varphi_{_{2}}\circ\varphi_{_{1}})(\mathbf{s}^{k-1})
\end{equation}
Starting with some admissible $\mathbf{s}^{0}$, we construct the sequence
$(\mathbf{s}^{k}\mathbf{)}$.

The convergence of the algorithm, taking into account Zangwill's global
convergence Theorem \cite{Zangwill 1969}, is justified in the same way as
the authors did in \cite{Bayon 2012}. The limit of
this ascending succession%
\begin{equation}
\underset{k\rightarrow\infty}{\lim}\mathbf{s}^{k}\label{conv}%
\end{equation}
provides the maximum.

\subsection{Step 2: The Unidimensional Problem}
Problem \eqref{UFunct}, \eqref{Ucon} is a unidimensional optimal
control problem with fixed end-time $T$, fixed initial state $x_{i}(0)$
and fixed end state $x_{i}(T)$, which can be expressed as
\begin{equation}%
\begin{array}
[c]{ll}
& \underset{u_{i}(t)}{\max}J_{i}^{\mathbf{s}}=\underset{u_{i}(t)}{\max}%
\int_{0}^{T}F_{i}^{\mathbf{s}}(t,x_{i}(t),u_{i}(t))dt\\
\text{subject to:} & \dot{x}_{i}(t)=f_{i}(t,x_{i}(t),u_{i}(t))\\
& x_{i}(0)=x_{i0};\text{ }x_{i}(T)=x_{iT}\\
  & u_{i}(t)\in\left[  h_{i\min},h_{i\max}\right]  \,\pedro{}
    \text{for each}\;t\in[0,T]\pedro{}
\end{array}
\label{Pr Pont}%
\end{equation}

where $f_{i}(t,x_{i}(t),u_{i}(t))=f_{l,i}(x_{i}%
)-g_{i}(s_{1},\cdots,s_{i-1},x_{i},s_{i+1},\cdots,s_{n})-u_{i}(t).$

Based on PMP \cite{Chiang 2000}, we establish the necessary
conditions of optimality for our unidimensional problem \eqref{Pr Pont}.
We require the following definition.

\textbf{Definition 2.} Let $x_{i}\in\Theta_{i}^{\mathbf{s}}$. The
\emph{coordination function} of $x_{i}$, ${\mathbb F}${}$_{x_{i}}(t)$ is the function in
$[0,T]$, defined as follows:%
\begin{equation}
{{\mathbb F}{}}_{x_{i}}(t)=\left(  F_{i}^{\mathbf{s}}(t,x_{i}(t),u_{i}%
(t))\right)  _{u_{i}}\cdot e^{\int_{0}^{t}\left(  f_{i}(t,x_{i}(v),u_{i}%
(v))\right)  _{_{x_{i}}}dv}+\int_{0}^{t}\left(  \left(  F_{i}^{\mathbf{s}%
}(t,x_{i}(v),u_{i}(v))\right)  _{_{x_{i}}}e^{\int_{0}^{v}\left(  f_{i}%
(t,x_{i}(z),u_{i}(z))\right)  _{_{x_{i}}}dz}\right)  dv\label{EQ COOR}%
\end{equation}
where sub-indices $(\cdots)_{u_i}$ or $(\cdots)_{x_i}$ mean partial
differentiation with respect to the corresponding coordinate.

\textbf{Theorem 1. A necessary maximum condition}

Let $u_{i}^{\ast}$ be the optimal control; let $x_{i}^{\ast}\in\widehat{C}^{1}$
be a solution of the above problem. Then there exists a constant $K\in\mathbb{R}$ such that:%
\begin{equation}%
\begin{tabular}
[c]{lll}%
\textit{If }$h_{i\min}<u_{i}^{\ast}<h_{i\max}$ & $\Longrightarrow$ &
${{\mathbb F}{}}_{x_{i}^{\ast}}(t)=K$\\
\textit{If }$u_{i}^{\ast}=h_{i\min}$ & $\Longrightarrow$ & ${{\mathbb F}{}
}_{x_{i}^{\ast}}(t)\leq K$\\
\textit{If }$u_{i}^{\ast}=h_{i\max}$ & $\Longrightarrow$ & ${{\mathbb F}{}
}_{x_{i}^{\ast}}(t)\geq K$%
\end{tabular}
\end{equation}
\textbf{Proof.}

Let $H$ be the Hamiltonian function associated with the problem
\begin{equation}
H(t,x_{i},u_{i},\lambda)=F_{i}^{\mathbf{s}}\left(  t,x_{i},u_{i}\right)
+\lambda\cdot f_{i}\left(  t,x_{i},u_{i}\right)
\end{equation}
where $\lambda(t)$ is the co-state variable.  \pedro{}Pontryagin's
Maximum Principle (PMP) asserts that\pedro{} in order for
$u_{i}^{\ast}\in\left[ h_{i\min},h_{i\max}\right] $ to be optimal, a nontrivial
function $\lambda$ must exist satisfying the following conditions:%
\begin{align}
\dot{x}_{i}^{\ast}  & =H_{\mathbf{\lambda}}=f_{i}\left(  t,x_{i}^{\ast}%
,u_{i}^{\ast}\right)  ;\text{  }x_{i}^{\ast}(0)=x_{i0}^{\ast},\text{  }%
x_{i}^{\ast}(T)=x_{iT}^{\ast}\label{eq:boundaries}\\
\dot{\lambda}  & =-\left(  H(t,x_{i}^{\ast},u_{i}^{\ast},\lambda)\right)
_{x_{i}} \label{lambda p}\\
H(t,x_{i}^{\ast},u_{i}^{\ast},\lambda)  & =\underset{u_{i}(t)\in\left[
h_{i\min},h_{i\max}\right]  }{\max}H(t,x_{i}^{\ast},u_{i},\lambda
)\label{Max H}%
\end{align}
\pedro{}Equations \eqref{eq:boundaries} and \eqref{lambda p} are a Boundary
Value Problem for $x^{\ast}_i$ and $\lambda$, which has solution as piecewise
$C^1$ functions. Specifically, $\lambda(t)$ satisfies:\pedro{}
\begin{equation}
\dot{\lambda}(t)=-\left(  H(t,x_{i}^{\ast},u_{i}^{\ast},\lambda)\right)
_{x_{i}}=-\left(  F_{i}^{\mathbf{s}}\right)  _{_{x_{i}}}-\lambda
(t)\cdot\left(  f_{i}\right)  _{_{x_{i}}}%
\end{equation}
and hence:
\begin{equation}
\lambda(t)=\left[  K-\int_{0}^{t}\left(  F_{i}^{\mathbf{s}}\right)  _{_{x_{i}%
}}e^{\int_{0}^{v}\left(  f_{i}\right)  _{_{x_{i}}}dz}dv\right]  e^{-\int
_{0}^{t}\left(  f_{i}\right)  _{_{x_{i}}}dv}\label{lambda}%
\end{equation}
where $K=\lambda(0)$.

From (\ref{Max H}) follows that, for each $t$ the value $u_{i}^{\ast}$
maximizes%
\begin{equation}
H(u_{i}):=H(t,x_{i}^{\ast},u_{i},\lambda)=F_{i}^{\mathbf{s}}\left(
t,x_{i}^{\ast},u_{i},u\right)  +\lambda\cdot f_{i}\left(  t,x_{i}^{\ast}%
,u_{i}\right)  ,\text{ }\forall u_{i}(t)\in\left[  h_{i\min},h_{i\max}\right]
\end{equation}

Bearing in mind that%
\begin{equation}
\frac{\partial H(u_{i})}{\partial u_{i}}=\left(  F_{i}^{\mathbf{s}%
}\right)  _{u_{i}}+\lambda(t)\left(  f_{i}\right)  _{u_{i}}\label{der}%
\end{equation}

there are three possibilities:

\begin{enumerate}
\item  $h_{i\min}<u_{i}^{\ast}<h_{i\max}$. In this case, $\frac{\partial
 H(u_{i})}{\partial u_{i}}=0.$ From (\ref{lambda}) and (\ref{der}) and
noticing that $\left(  f_{i}\right)  _{_{u_{i}}}=-1$, we obtain:%
\begin{equation}
0={{\mathbb F}{}}_{x_{i}^{\ast}}(t)-K\Rightarrow{{\mathbb F}{}}%
_{x_{i}^{\ast}}(t)=K\label{T1}%
\end{equation}

\item $h_{i\min}=u_{i}^{\ast}$. In this case, $\frac{\partial H(u_{i}%
)}{\partial u_{i}}\leq0.$ an analogous reasoning gives:
\begin{equation}
{{\mathbb F}{}}_{x_{i}^{\ast}}(t)\leq K\label{T2}%
\end{equation}

\item $u_{i}^{\ast}=h_{i\max}$. In this case, $\frac{\partial H(u_{i}%
)}{\partial u_{i}}\geq0$ and the same argument provides:
\begin{equation}
{{\mathbb F}{}}_{x_{i}^{\ast}}(t)\geq K\label{T3}%
\end{equation}
\end{enumerate}
This ends the proof.
$\hfill\square$

From the computational point of view, the construction of the solution
consists of two main steps which we proceed to explain.

\subsection{Step 3: The construction of $x_{i}^{K}$}

In this section, we \pedro{}describe how, using the $i$-th maximizing map, we compute an approximation to the function $x_i^K$ for a given value $K$,
using the fact that
\pedro{}
for each $i$, $x_{i}^{K}$ maximizes Problem (\ref{Pr Pont}). This \emph{approximate} construction of $x_{i}^{K}$ can be performed using a
discretized version of the following equation, which we shall call
the \emph{coordination equation}:%
\begin{equation}\label{EQ-COOR-K}
K=\left(  F_{i}^{\mathbf{s}}\right)  _{u_{i}}\cdot e^{\int_{0}^{t}\left(
    \pedro{}f_{i}\pedro{}\right)  _{_{x_{i}}}dv}
+\int_{0}^{t}\left(  F_{i}^{\mathbf{s}}\right)
_{_{x_{i}}}e^{\int_{0}^{v}\left(  \pedro{}f_{i}\pedro{}\right)  _{_{x_{i}}}dz}dv
\end{equation}
For given $K$, we compute $x_{i}^{K}$, \pedro{}by means of (\ref{EQ-COOR-K})\pedro{}. When the
values obtained do not obey the control constraints, we force the
solution to belong to the boundary until the moment established by conditions
(\ref{T2}) and (\ref{T3}). This computation is carried out using
polygonals (an adaptation of Euler's method).

In summary, we divide the interval $[0,T]$ using $N$ nodes:%
\begin{equation}
\pedro{}t_{0}<t_{1}<\dots<t_{N-1}\pedro{}
\end{equation}
with $n=0,..., N-1$ so that $t_{0}=0$ and $t_{N-1}=T$.

Next, we solve the coordination equation \eqref{EQ-COOR-K} at each node.
Starting at $t_{0}=0$, where $x_{i0}$ is known, we compute, from
the discretized version of the coordination equation \eqref{EQ-COOR-K},
the corresponding control $u_{i0}$. Using this value for the control, we calculate the derivative of the state
equation
\begin{equation}
  \dot{x}_{i}(t)= \pedro{}
  f_i(t,x_i(t),u_i(t)) =
  f_{l,i}(x_{i})-g_{i}(s_{1},\dots,s_{i-1},x_{i},s_{i+1},\dots,s_{n})-u_{i}(t)
\end{equation}
In this case:%
\begin{equation}
  \label{eq:numerical-diffeq}
\dot{x}_{i0}=\pedro{}f_{i}\pedro{}\left(  t,x_{i0},u_{i0}\right)
\end{equation}
And using this value we jump to the next node, using Euler's method, imposing
the condition:
\begin{equation}
x_{i1}=x_{i0}+\pedro{}f_{i}\pedro{}\left(  t,x_{i0},u_{i0}\right)   d
\end{equation}
And, in general:
\begin{equation}
  x_{ik}=x_{ik-1}+\pedro{}f_i\pedro{}\left(  t,x_{ik-1},u_{ik-1}\right) d
\end{equation}
where the step $d$ is:
\begin{equation}
d=\frac{T}{N}%
\end{equation}
%

%
That way,\pedro{} for $k=0,\dots, N-1$, we go through\pedro{} all the nodes until the whole interval $[0,T]$ is covered, linking the arcs continuously.

We must also take into account the restrictions on the control:
\begin{equation}
  h_{i\min} \leq u_i(t) \leq h_{i\max}  ,\text{ }0\leq t\leq T.
\end{equation}
\pedro{}When these are not verified by the value obtained from solving the coordination equation \eqref{EQ-COOR-K}, \pedro{}we force the solution to remain on the border until the moment established by conditions  \eqref{T2} and \eqref{T3}:%
\begin{equation}
\mathbb{F}_{x_{i}}(t)\leq K
\end{equation}%
\begin{equation}
\mathbb{F}_{x_{i}}(t)\geq K
\end{equation}
\bigskip

\subsection{Step 4: Computation of the optimal $K$}

The computation of the optimal $K$ is achieved by means of an adaptation
of the shooting method. We search for the extremal $x_i^K$ fulfilling the second boundary condition by varying the coordination constant $K$
(\ref{CT}). The procedure is similar to the shooting method used to solve
second-order differential equations with boundary conditions, which may be
performed approximately using elemental procedures. Starting from two
values for the coordination constant, $K$: $K_{\min}$ and $K_{\max}$ and using
a conventional method such as the secant method, our algorithm converges
satisfactorily, as we shall see in the examples.

In summary, the problem consists in finding, for each $K$, the function $x_{i}^{K}$ satisfying the initial condition
\begin{equation}
x_{i}^{K}(0)=x_{i0}%
\end{equation}
the conditions of Theorem 1 and, among these, the one that satisfies the
second boundary condition:
\begin{equation}
x_{i}^{K}(T)=x_{iT}\label{CT}%
\end{equation}

Consider the \emph{shooting function}:
\begin{equation}
\varphi(K):=x_{i}^{K}(T)
\end{equation}
Solving our final condition problem is equivalent to solving the (finite dimensional) equation

\begin{equation}
\varphi(K)-x_{iT}=0
\end{equation}
We solve this equation by means of the secant method, imposing a stopping
criterion based on a pre-established tolerance:
\begin{equation}
error=|\varphi(K)-x_{iT}|\leq tol_{K} \label{Tol k}%
\end{equation}
The secant method gives:
\begin{equation}
  \label{eq:secant}
K_{j+1}=K_{j}-\frac{K_{j}-K_{j-1}}{\varphi(K_{j})-\varphi(K_{j-1})}%
\varphi(K_{j})
\end{equation}%

Thus, for each $i$-th component, we construct a sequence
$\{K_{j}\}_{j\in\mathbb{N}}$ such that $x_{i}^{K_{j}}(T)$ converges to
$x_{iT}$. The optimal $K$ is the one verifying
\begin{equation}
\underset{j\rightarrow\infty}{\lim}K_{j}=K
\end{equation}
with $x_{i}^{K}(T)=x_{iT}$.

The stopping criterion for the algorithm is based on the desired
tolerance. We shall see this fact in greater detail
in Example 1 in Section 4.

\begin{figure}[h!]
  \centering
  \includegraphics{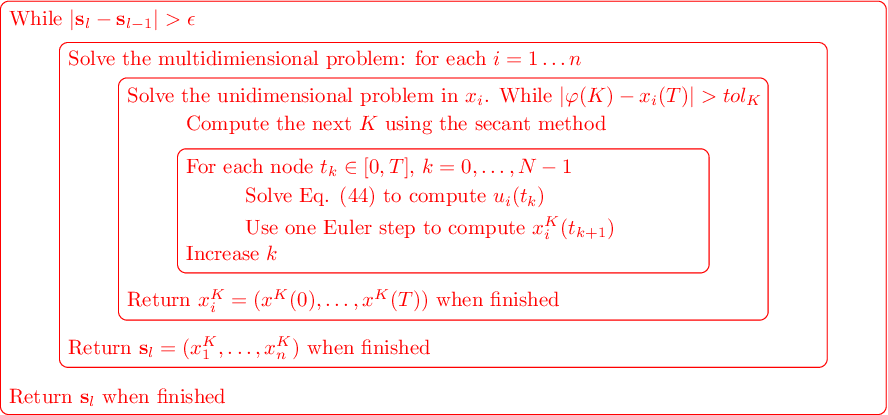}
  \caption{\color{black}Schematic flow of the whole algorithm.\color{black}}
  \label{fig:flow-chart}
\end{figure}

Figure \ref{fig:flow-chart} shows a schematic picture of our algorithm.

\color{black}

\section{Numerical Examples}
We now show some numerical simulations, with a twofold aim: firstly,
to illustrate how our algorithm works; secondly, to analyze the
influence on the solution of some of the parameters of the problem
and/or of the model.
\subsection{Example 1: Base Case}

Let us begin with what we shall call the \emph{base case} and which will
allow us to compare the results obtained with different parameters. This
is the $3$-species model described in \eqref{SD1}, \eqref{SD2} and
\eqref{SD3}. The values of the cost model are summarized in Table I,
whereas Table II shows the parameters of the biological model, including
the interactions among the species. Growth is measured in $10^6\mathrm{kg}$,
and time in units of years. Prices of the economic parameters are in \euro/kg
and costs are in $10^6$\euro.

\begin{table}[h!]
  \centering
\begin{tabular}
[c]{c}%
Table I. Parameters of the cost model.\\%
\begin{tabular}
[c]{ccccc}%
Specie & $p_{0}$ & $p_{1}$ & $c$ & $\alpha$\\\hline
1 & $0.9$ & $0.01$ & $75$ & $1.1$\\
2 & $1.9$ & $0.02$ & $85$ & $1.2$\\
3 & $2.8$ & $0.03$ & $60$ & $1.4$%
\end{tabular}
\end{tabular}
\end{table}

\begin{table}[h!]
  \centering
\begin{tabular}
[c]{c}%
Table II. Parameters of the biological models.\\%
\begin{tabular}
[c]{c|lllll}%
$f_{l,i}(x_{i}(t))$ &  &  & $g_{i}(\mathbf{x(}t\mathbf{)})$ &  & \\\hline
\multicolumn{1}{l|}{$r_{1}=0.5$} & $c_{12}=2.10^{-4}$ & $c_{21}=10^{-5}$ &
$c_{31}=10^{-4}$ & $\gamma_{12}=1.1$ & $\gamma_{123}=1.1$\\
\multicolumn{1}{l|}{$r_{2}=0.3$} & $c_{13}=3.10^{-5}$ & $c_{23}=2.10^{-5}$ &
$c_{32}=10^{-4}$ & $\gamma_{13}=1.0$ & $\gamma_{132}=1.2$\\
\multicolumn{1}{l|}{$r_{3}=0.2$} & $c_{123}=10^{-8}$ & $c_{213}=10^{-7}$ &
$c_{312}=0$ & $\gamma_{21}=1.0$ & $\gamma_{213}=1.0$\\
\multicolumn{1}{l|}{$k_{1}=1000$} & $\beta_{12}=1.0$ & $\beta_{21}=1.0$ &
$\beta_{31}=1.0$ & $\gamma_{23}=1.2$ & $\gamma_{231}=1.0$\\
\multicolumn{1}{l|}{$k_{2}=700$} & $\beta_{13}=1.2$ & $\beta_{23}=1.0$ &
$\beta_{32}=1.0$ & $\gamma_{31}=1.1$ & $\gamma_{312}=0.$\\
\multicolumn{1}{l|}{$k_{3}=600$} & $\beta_{123}=0.9$ & $\beta_{213}=1.0$ &
$\beta_{312}=0.$ & $\gamma_{32}=1.0$ & $\gamma_{321}=0.$%
\end{tabular}
\end{tabular}
\end{table}

The values in Table I are not taken from any real species but they are
inspired by \cite{Agnarsson 2008}, where a detailed modeling was carried
out for cod, capelin and herring in Norway, Iceland and Denmark. We
have respected the ratio between income and cost in the revenue function
which they checked against the real revenues. As for the parameters in
Table II about species interaction, $g_{i}(\mathbf{x(}t\mathbf{)})$
and the simple logistic function $f_{l,i}(x_{i}(t))$, the choice of parameters
follows the ideas in \cite{Agmour 2017}. All our results were obtained
using a custom-written software implemented in Mathematica
10.0\copyright.
\begin{figure}[h!]
  \centering
{\includegraphics[
height=1.5298in,
width=2.5201in
]%
{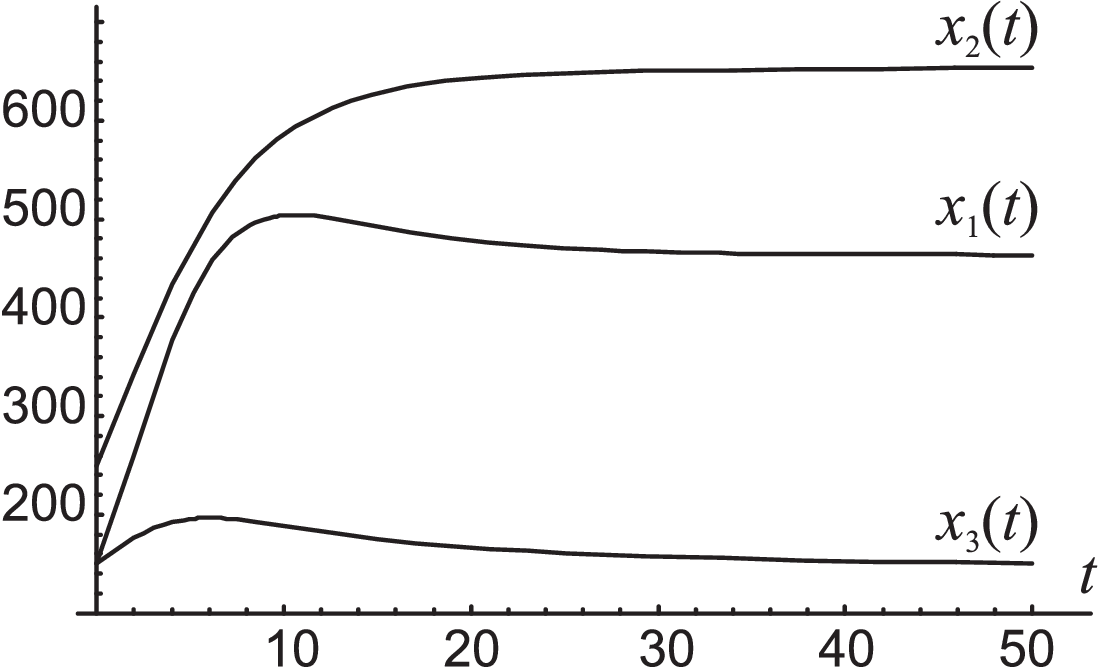}%
}
\\
\text{Fig. 2. Steady state of biological equilibrium, no harvesting.}%
\end{figure}
Figure 2 shows the solution obtained in the absence of harvesting, that is,
the long-term steady-state solution of the biological model \eqref{SD1},
\eqref{SD2}, \eqref{SD3}, with no optimization of the functional. Starting
from the initial stock values (in $10^{6}$kg):
\begin{equation}
x_{1}(0)=150;\text{ \ }x_{2}(0)=250;\text{ \ }x_{3}(0)=150\label{X0}
\end{equation}
the values obtained for the steady state of the stock ($10^6$kg) of
each species after $100$ years are:
\begin{equation}
x_{1}^{\ast}=463.69;\text{ \ }x_{2}^{\ast}=654.65;\text{ \ }x_{3}^{\ast
}=146.96\label{STS}%
\end{equation}
However, our interest is not the steady-state: starting from the known
initial stocks \eqref{X0}, we fix a time horizon of say $T=10$ years and,
for different reasons, (environmental, economical, ecological,
legislative...), we state a final stock for each species to reach in
that time-span, different from the one set by the free biological
evolution of the system \eqref{STS}. Let us set, in this example,
the following target values:
\begin{equation}
x_{1}(T)=500;\text{ \ }x_{2}(T)=500;\text{ \ }x_{3}(T)=200
\end{equation}
The solution will optimize the harvesting revenue during that time,
by choosing the optimal controls, i.e. how each species should be
fished. We set a $\delta=5$ percent discount rate.

The results obtained for the optimal stock profile, $x(t)$ (in
$10^{6}$kg), and the optimal harvest path, $h(t)$
(in $10^{6}$kg/year), of each species are shown in Figs. 3 and 4
respectively.
\begin{figure}[h!]
  \centering
{\includegraphics[
height=1.2332in,
width=6.3123in
]%
{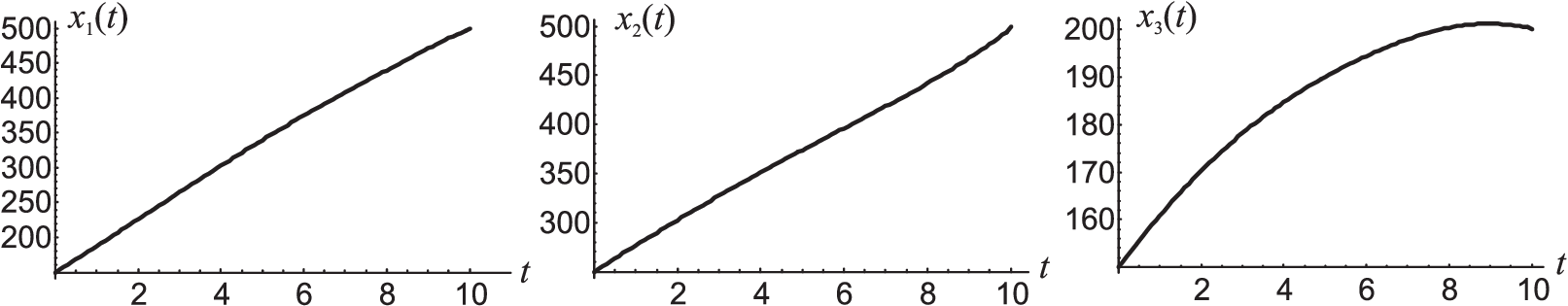}%
}
\\
\text{Fig. 3. Optimal solution for the stock profile.}%
\end{figure}%
\begin{figure}[h!]
  \centering
{\includegraphics[
height=1.2496in,
width=6.0416in
]%
{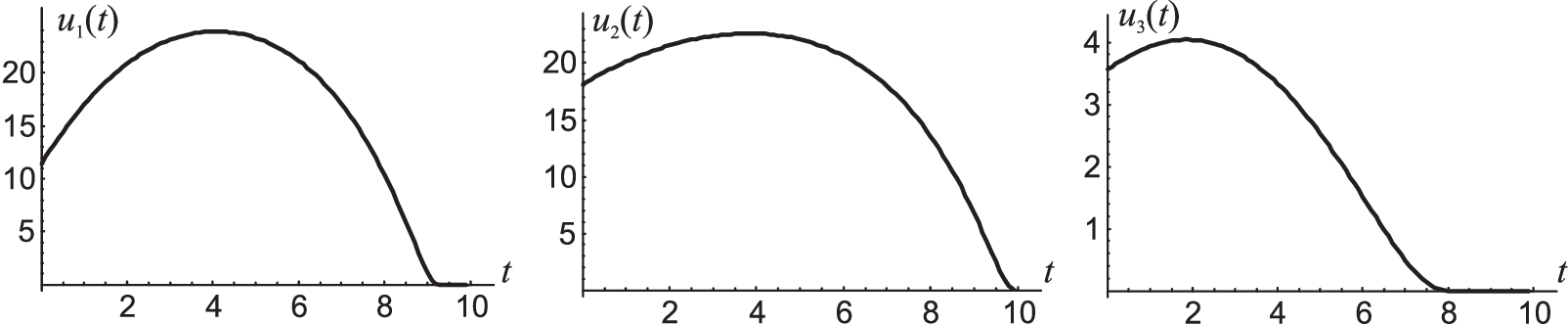}%
}
\\
\text{Fig. 4. Optimal solution for the harvest profile.}%
\end{figure}
In computing the solution, we used a discretization of $100$ sub-intervals.
For each species and for the secant method, we chose a tolerance \eqref{Tol k}
$tol_K = 5\cdot 10^{-2}$. The stopping criterion for the multi-dimensional
algorithm was that the difference between two consecutive iterations of
the values of $K$ was less than $10^{-6}$. With these parameters, the
algorithm converges in just $9$ iterations to the final solution (see
Fig. 5).
\begin{figure}[h!]
  \centering
{\includegraphics[
height=1.4581in,
width=2.7544in
]%
{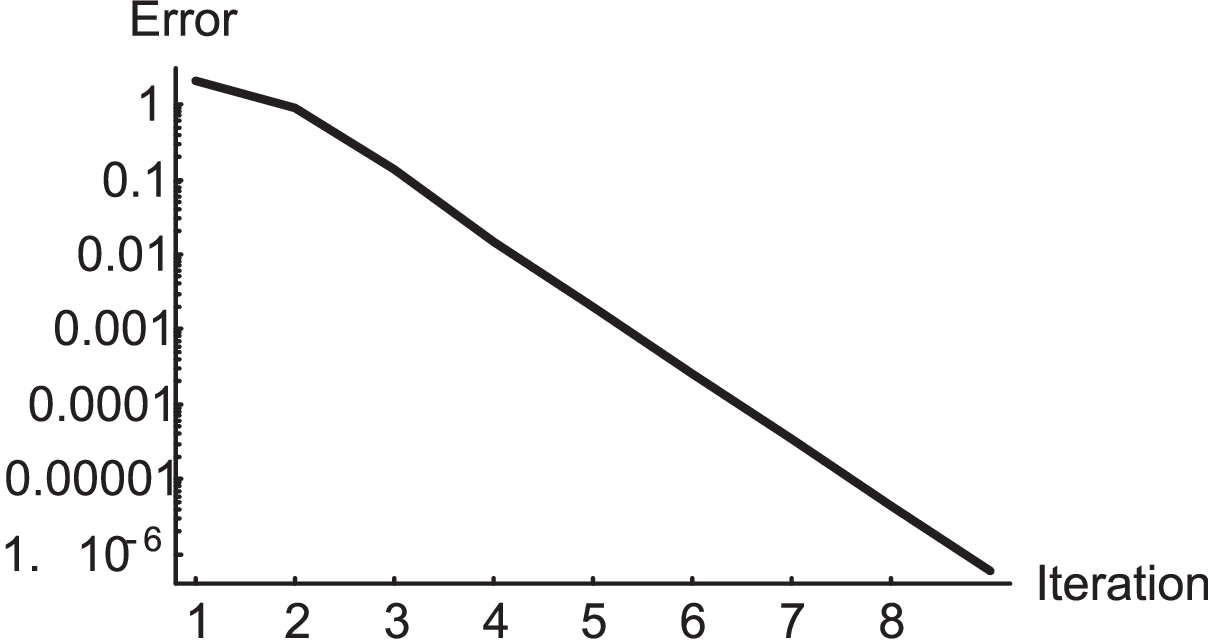}%
}
\\
\text{Fig. 5. Convergence of the multi-dimensional algorithm.}%
\end{figure}
The optimal values of $K$ for each species are, respectively:
\begin{equation}
  K_{1}^{\ast}=-0.0299251;\text{ \ }K_{2}^{\ast}=0.4474817;
  \text{ \ }K_{3}^{\ast}=1.6556858
\end{equation}
and the final revenue after $10$ years is  $235.381$(in $10^{6}$\euro).
In finding the solution, we have considered the control limited by
$u_{i}(t)\in\left[  0,25\right]$ (in $10^6$kg). Fig. 4 shows how the
maximum is not reached at any moment but the minimum is. For species
1, harvesting stops in the interval $[9.2,10]$, and for species 3,
in the interval $[8,10]$. We consider these two facts quite remarkable,
as they state roughly that the optimum is reached by fishing, mainly,
in the middle part of the optimization interval and practically stopping the
harvest at its end. At this point, and almost without the intervention
of the control, the system reaches the desired final state ``by itself''.

\subsection{Example 2: Influence of final stock and T}

In this section we fix the biological parameters for each species, which
is a fairly realistic assumption for our mid-term horizon. For the same
reason, we fix the economic and price parameters. In this second example,
we are going to study a problem we deem relevant: given some initial
stocks, what are the possible achievable final values for them and
what influence does the time-span have on them?

\paragraph{Varying the coordination constant $K$} Let us, for the sake of concreteness, assume that the final stocks of
species $2$ and $3$ are imposed and are the same as in the previous
example: $x_{2}(T)=500$ and $x_{3}(T)=200$ (in $10^{6}$kg).
For species $1$, starting from the initial stock $x_{1}(0)=150$ (in $10^6$kg),
what are the final stock values $x_1(T)$ that we can impose that
are biologically possible after $T=10$ years? For the sake of simplicity,
we shall assume at this point that there are no
upper limits for the controls.

In order to answer the question above, we are going to perform an
adaptation of the algorithm previously presented. We are going to
consider the one-dimensional problem on species $1$, imposing, as
known stock values, for each $t$, those obtained in the base example
above, $x_2(t)$ and $x_3(t)$.

Then, instead of fixing the final state $x_1(T)$, we are going to
carry out a sweeping of the admissible interval for the coordination
constant $K$, using the shooting function:
\begin{equation}
\varphi(K):=x_{1}^{K}(T)
\end{equation}
Table III shows the results obtained for the final stock of $x_1(T)$ in
$10^6$kg when sweeping the $K$-interval $[-2,2]$.
\begin{table}[h!]
  \centering
\begin{tabular}
[c]{c}%
Table III. Sweeping of the $K$-values in the shooting function $\varphi(K)$.\\
\begin{tabular}
[c]{c|ccccccccccc}%
$K$ & $-2$ & $-1.6$ & $-1.2$ & $-0.8$ & $-0.4$ & $0$ & $0.4$ & $0.8$ & $1.2$ &
$1.6$ & $2.0$\\\hline
$x_{1}^{K}(T)$ & 0.03 & 0.06 & 0.16 & 1.16 & 141.26 & 340.57 & 521.67 &
578.65 & 580.20 & 580.20 & 580.20
\end{tabular}
\end{tabular}
\end{table}
Notice how the admissible values for $x_1^K(T)$ is the interval
$[0.03,580.20]$ (in $10^6$kg). Fig. 6 shows the optimal profiles of
$x_1^K(t)$ for the values of $K$ given in Table III, whereas
Fig. 7 shows the optimal controls obtained for those values.
\begin{gather*}
{\includegraphics[
height=2.047in,
width=5.4656in
]%
{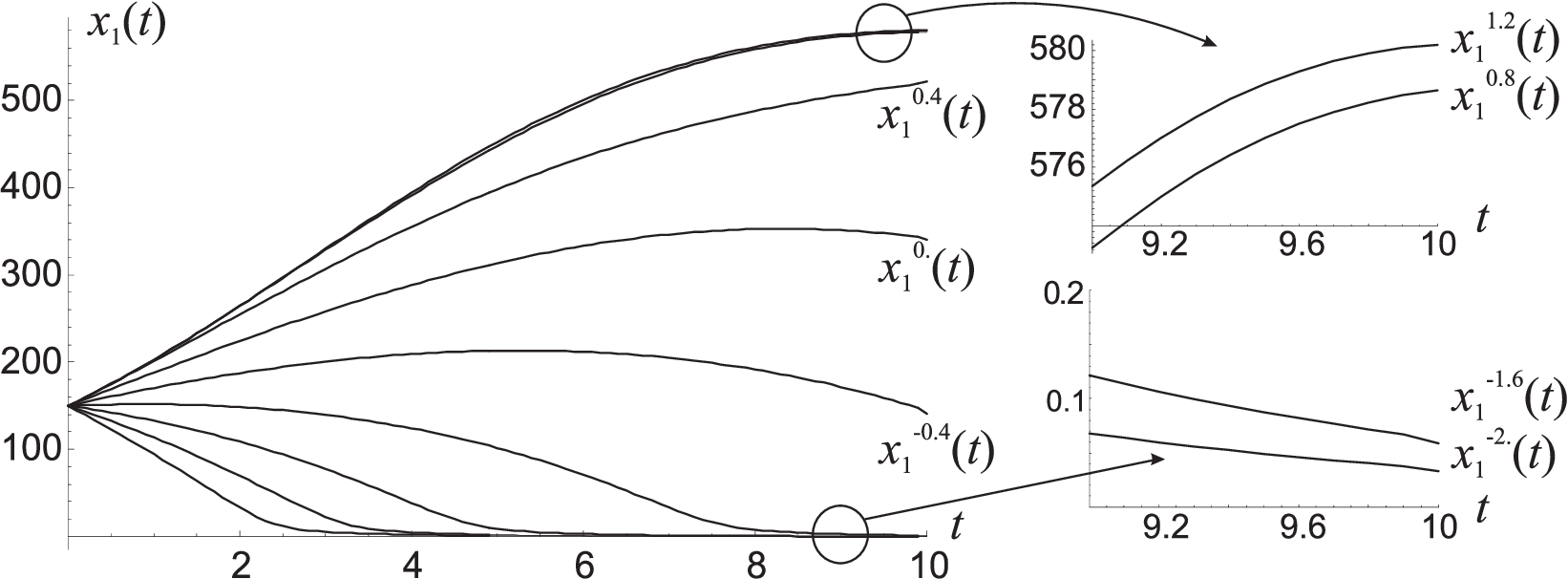}%
}
\\
\text{Fig. 6. Admissible interval for the values of the final stock.}%
\end{gather*}
As one can see, we have constructed the extremal field for the
variable $x_1(t)$: the zooms on the right side shows how the extremals
do not meet.

The maximum value for the stock is quickly reached for values of $K$
greater than $1.2$ and it corresponds to those cases where the harvesting
(the control) is $0$ (see Fig. 7). On the other hand, the minimum
value of the stock (essentially $0$) is reached with the optimal
harvesting profiles shown in Fig. 7. At the beginning of the
interval, the fishing activity is remarkably large, so that when the
stock of species $1$ has already been greatly diminished, only the
competency effects of species $2$ and $3$ and practically no
harvesting, species $1$ gets to its practical biological demise.
Obviously, this is not the desired value at all: we show it as a
validation of the theoretical soundness of our algorithm.
\begin{gather*}%
{\includegraphics[
height=1.9744in,
width=3.3295in
]%
{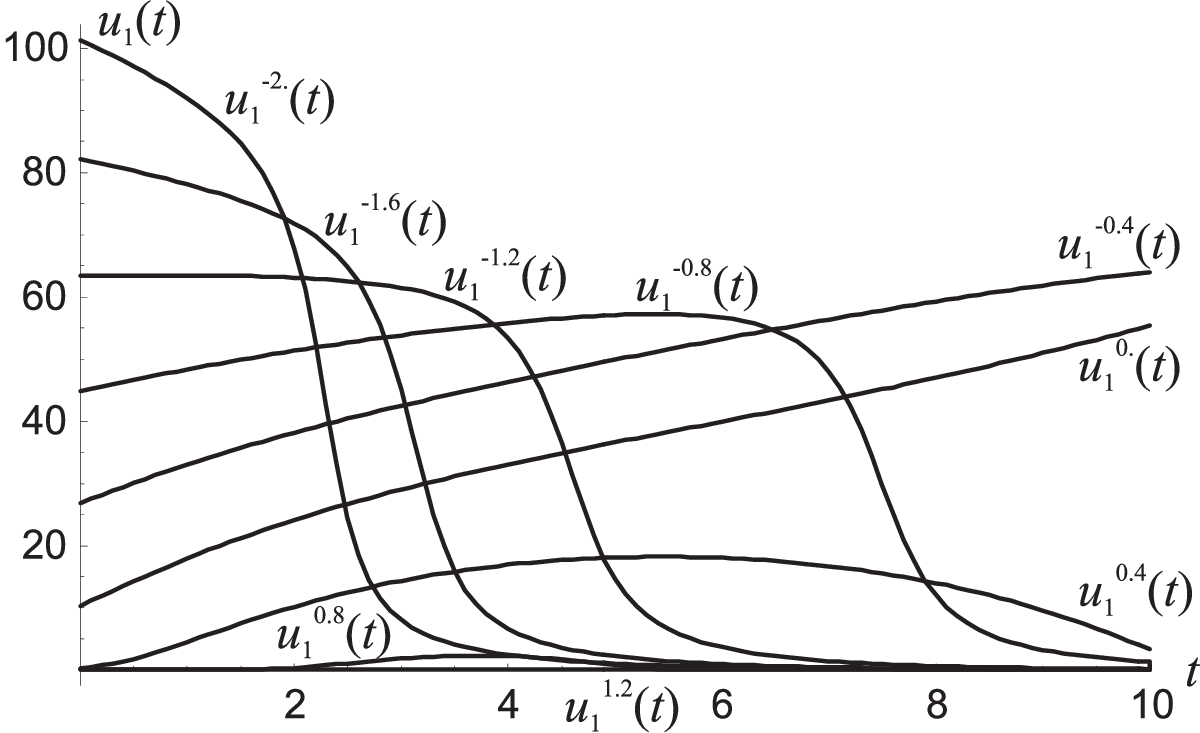}%
}
\\
\text{Fig. 7. Optimal controls for the sweeping of the $K$-values.}%
\end{gather*}

\paragraph{Influence of $T$ on the admissible interval} This is the
second question we are going to address. Table IV summarizes the results
for the values $T=8, 10, 12$ and $14$ years: the maximum achievable value
of $x_1(T)$ for each $T$ is shown (in $10^6$kg), corresponding to the
maximum values of $K$. Notice that the minimum achievable value of
$x_1(T)$ is the same in each case, essentially $0$ as shown in the previous
example.
\begin{table}[h!]
  \centering
\begin{tabular}
[c]{c}%
\text{Table IV. Influence of $T$ on the admissible interval.}\\%
\begin{tabular}
[c]{ccccc}%
$T$ & $8$ & $10$ & $12$ & $14$\\\hline
$\max x_{1}^{K}(T)$ & 538.09 & 580.20 & 598.88 & 605.42
\end{tabular}
\end{tabular}
\end{table}
Notice how, the longer the time $T$, the larger the admissible variation
rank for the final stock of species $1$, $x_1(T)$ which is reachable. However,
Fig. 8 shows clearly how these increase is not linear. Actually, the species
has a biological limit and, notwithstanding the time-span, the maximum admissible
stock does not grow indefinitely. We consider this result (which is
reasonable from the biological point of view) provides an important
contingency of which the harvesting managers have to be aware for the
mid-term planning.
\begin{gather*}%
  {\includegraphics[
height=1.3in
]%
{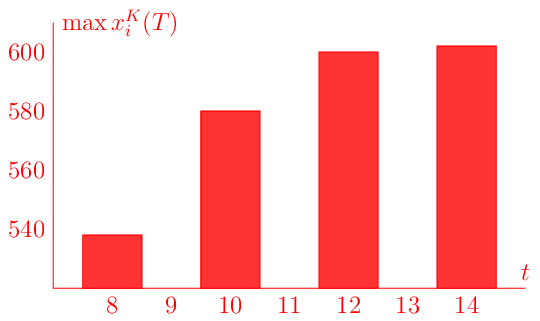}%
}
\\
\text{\pedro{}Fig. 8. Histogram corresponding to Table IV: influence of $T$ on $\max x_1^K(T)$.\pedro{}}%
\end{gather*}

We finish this section by remarking that what we have obtained is
just an \emph{estimation} of the structure of the set of admissible
solutions for $x_1(T)$, as the values of $x_2(t)$ and $x_3(t)$ were
chosen from the base case. These values, even starting from the
same initial stock value (each its own) and reaching the same
final value (ibid.), might take different intermediate values and
this would, certainly, influence the admissible interval of $x_1(T)$.
However, several numerical simulations lead us to conjecture that
this variation should be quite small and that our estimation is
rather good.

\section{Conclusions and future perspectives}

Going further than the study of the long-term steady-state solution of
the multi-species problem in the renewable resources problem, we have
considered a joint finite mid-term horizon and a stock target at the
end of the optimization interval. The solution of this problem involves
a highly dynamic system which requires the development of a method
allowing to solve it satisfactorily. This algorithm, an adaptation of
the cyclic coordinate descent provides a procedure whose
complexity is not essentially influenced by the dimensionality
of the problem, as it reduces to an iteration of one-dimensional
algorithms. Each of these is tackled by approximately solving what we
call the \emph{coordination equation}. This, together with a
customization of the shooting method, eliminates the multi-dimensionality
and thus, removes complexities due to the biological interactions
from the resolution process.

\luis{}
In brief, we consider that the main novel characteristics of our work are:
\begin{itemize}
\item The setting of a fixed stock at a mid-term horizon as one of the contraints.
\item The arbitrary dimensionality of the considered problem: we allow for any number of species, with arbitrary simultaneous interactions.
\item The flexibility of the model, allowing for nonlinearity in the interactions, which greatly increases the modeling possibilities.
\item The method used to solve the problem: despite the nonlinearity of the model, we have an algorithm whose complexity is not inherently increased by either the dimensionality or the interactions.
\item Our model permits studying problems of very different nature: predator-prey models, competition, host-parasitoid, infectious disease, etc.
\end{itemize}
\luis{}

As for the solution we obtain, we remark that it is possible to achieve
the desired final stock values while maximizing the revenues of harvesting
during the optimization interval. As a matter of fact, we observe that
the optimum is reached in such a way that the harvesting practically stops
near the end of the interval and the biological system evolves by itself
to the desired outcome. We have also analyzed the variation of the
range of admissible values for the stock depending on the duration
of the optimization interval: this gives the policy makers an a priori
set of limit values of each species which they can reach without
harm to the biological system. One remaining question which we do not
know how to tackle yet is the structure of the joint set of admissible
values for each species at the same time. This is an open question which
we think requires a richer set of tools or a deeper insight than we
have at present.

\end{document}